\documentclass[10pt,letterpaper,reqno]{amsart}

\usepackage[left=1.25in, right=1.25in, top=1.25in, bottom=1.5in]{geometry}
\usepackage{amsthm,amsmath,amssymb,amsfonts}
\usepackage{hyperref}
 \hypersetup{colorlinks=true,linktocpage=true,linkcolor=blue,citecolor=blue,urlcolor=blue}
\usepackage[utf8x]{inputenc}
\usepackage[all]{xy}
\usepackage{enumerate}
\usepackage{lipsum}
\usepackage{stmaryrd}
\usepackage{pdfpages}
\usepackage{mathrsfs}

\hypersetup{
  pdftitle={The Woods Hole trace formula and indices for vector fields},
  pdfauthor={Valente Ramirez}
}

\newcommand{\C}{\mathbb{C}}
\renewcommand{\P}{\mathbb{P}}

\renewcommand{\O}{\mathcal{O}}
\newcommand{\F}{\mathcal{F}}
\newcommand{\Fs}{\mathscr{F}}
\newcommand{\tr}{\operatorname{tr}}

\newtheorem{theorem}{Theorem}
\newtheorem{corollary}{Corollary}
\newtheorem{lemma}{Lemma}
\newtheorem{proposition}{Proposition}

\newtheorem{application}{Application}
  \theoremstyle{definition}
\newtheorem{definition}{Definition}

  \theoremstyle{remark}
\newtheorem{remark}{Remark}

\title[The Woods Hole trace formula and indices for vector fields]{The Woods Hole trace formula and indices for\\ vector fields and foliations on $\C^2$}
\author[Valente Ram\'{i}rez]{Valente Ram\'{i}rez}
\address{Department of Mathematics, Cornell University\\ 108 Malott Hall\\ Ithaca, NY\\ 14850}
\email{\href{mailto:valente@math.cornell.edu}{valente@math.cornell.edu}}
\date{}
\thanks{This work was supported by the grants UNAM-DGAPA-PAPIIT IN-102413 and CONACYT 219722.}

\begin{document} 

\begin{abstract}
The purpose of this work is to introduce the reader to the \emph{Woods Hole trace formula} and to show how several well-known index theorems for foliations, vector fields and endomorphisms can be rephrased as particular cases of such formula. The main objective in doing this is to use the Woods Hole formula in the future to produce new index theorems.
\end{abstract}

\maketitle

\vspace{-6mm}
\setcounter{section}{-1}
\section{Motivation}

Let $v$ be a polynomial vector field on $\C^2$ having only non-degenerate singularities. Assume $v$ is \emph{non-dicritic}, i.e.~the line at infinity $L$ is invariant for the induced foliation $\F_v$ on $\P^2$ and assume also that the singularities of $\F_v$ along $L$ are non degenerate -- note that all these assumptions are generic conditions for polynomial vector fields on $\C^2$. Each singular point $p$ of $v$ carries two numeric invariants, the eigenvalues of the linearization matrix $Dv(p)$. On the other hand, each singular point at infinity carries only one invariant, its Camacho-Sad index. We know that there are certain relations among all these invariants, namely the \emph{Euler-Jacobi relations}:
\begin{eqnarray} 
 \sum_{v(p)=0}\frac{1}{\det(Dv(p))}=0, \label{eq:EJ1} \\
 \sum_{v(p)=0}\frac{\tr(Dv(p))}{\det(Dv(p))}=0, \label{eq:EJ2}
\end{eqnarray}
the Baum-Bott formula:
\begin{equation} \label{eq:BB}
 \sum_{\operatorname{Sing }\F_v}\operatorname{BB}(\F_v,p)=(d+2)^2,
\end{equation}
and the Camacho-Sad formula for the invariant line $L$:
\begin{equation}
\sum_{L\cap\operatorname{Sing }\F_v}\operatorname{CS}(\F_v,L,p)=1. 
\end{equation}

The above relations are well-known, yet a simple dimension count shows that there must be more, currently unknown, such relations. The objective of this paper is to show, in \hyperref[sec:foliations]{\S\ref*{sec:foliations}} and \hyperref[sec:vectorfields]{\S\ref*{sec:vectorfields}}, how all the above relations can be reduced to a particular case of the Woods Hole formula. Our hope is to use the Woods Hole formula in the future to discover our missing relations. Some of these motivating questions are discussed in \cite{TwinVectorFields} and \cite{Guillot2004}.

In order to deduce the above index-theorems from Woods Hole we will first introduce the formula in \hyperref[sec:WoodsHole]{\S\ref*{sec:WoodsHole}}, and use use it to get index theorems for regular endomorphisms of $\P^2$ in \hyperref[sec:endomorphisms]{\S\ref*{sec:endomorphisms}}.

\section{The Woods Hole formula}\label{sec:WoodsHole}

The so-called \emph{Woods Hole formula} (also known as Atiyah-Bott fixed point theorem) is a generalization of the Lefschetz fixed point formula due to Atiyah and Bott in the complex analytic case \cite{AtiyahBott1966,AtiyahBott1967,AtiyahBott1968}, and to Verdier in the algebraic case \cite{Verdier1967} (see also \cite{Taelman2015,Beauville1972}). 

\medskip
Let $X$ be a compact complex manifold and $f$ a holomorphic endomorphism of $X$. Let $\O_X$ denote the structure sheaf on $X$ and let $\Fs$ be a  coherent analytic sheaf.

\begin{definition}\label{def:lift}
A \emph{lift} of $f$ to $\Fs$ is a morphism of $\O_X$-modules
\[ \varphi\colon f^*\Fs \to \Fs. \]
\end{definition}

Note that a lift of $f$ to $\Fs$ induces linear maps on cohomology $\varphi^{k}\colon H^{k}(X,f^*\Fs) \to H^{k}(X,\Fs)$. In addition, the pullback functor $f^*$ gives maps $H^{k}(X,\Fs) \to H^{k}(X,f^*\Fs)$. Composition of these maps gives endomorphisms on the cohomology groups:
\[ \xymatrix{
 \widetilde{\varphi}^k\colon H^{k}(X,\Fs) \ar[r]^{f^*} & H^{k}(X,f^*\Fs) \ar[r]^{\varphi^k} & H^{k}(X,\Fs)}.
\]

\begin{definition}\label{def:Lefschetznumber}
The \emph{Lefschetz number} $\mathrm{L}(f,\varphi,\Fs)$ of a lift $\varphi$ is defined by
\[ \mathrm{L}(f,\varphi,\Fs) = \sum_k (-1)^k\, \tr(\widetilde{\varphi}^k,\,H^{k}(X,\Fs)). \]
\end{definition}

Let $x\in X$. We denote by $\Fs_x$ the stalk of $\Fs$ at $x$, and by $\Fs(x)$ the fiber of $\Fs$ at $x$, namely
\[ \Fs(x)=\Fs_x\otimes_{\O_{X,x}}(\O_{X,x}/\mathfrak{m}_x). \]
Note that each fiber is a finite dimensional complex vector space. Since $\varphi$ induces morphisms on stalks
\[ \varphi_x\colon (f^*\Fs)_x \cong \Fs_{f(x)} \longrightarrow \Fs_x, \]
we see that if $p$ is a fixed point of $f$ we obtain a morphism $\varphi_p\colon\Fs_p\to\Fs_p$. In a similar way we obtain linear endomorphisms on the fibers $\varphi(p)=\Fs(p)\to\Fs(p)$. Finite dimensionality of the fibers guarantees that these endomorphisms have a well defined trace.

\begin{definition}\label{def:transversal}
We say that an endomorphism $f$ of $X$ is \emph{transversal} if its graph is transversal to the diagonal in $X\times X$. Equivalently, $f$ is transversal if and only if it has isolated simple fixed points only.
\end{definition}

\begin{theorem}[Woods Hole formula]\label{thm:WoodsHole}
If $f$ is a transversal endomorphism then
\[ \sum_{f(p)=p} \frac{\tr(\varphi(p),\,\Fs(p))}{\det(1-df(p),\,\Omega^1_X(p))} = \mathrm{L}(f,\varphi,\Fs), \]
where $df\colon f^*\Omega^1_X\to\Omega^1_X$ denotes the differential of $f$.
\end{theorem}

The Woods Hole formula is commonly used in the particular cases where $\Fs=\Omega^k_X$ and $\varphi=\wedge^k df$ (being the case $k=0$ the classical holomorphic Lefschetz formula). In such case, if $J_f$ denotes the Jacobian matrix of $f$, the formula becomes
\[ \sum_{f(p)=p}\frac{\tr(\wedge^k J_f(p))}{\det(I-J_f(p))} = \operatorname{L}(f,\wedge^k df,\Omega^k_X). \]
Note that the indices on the left hand side depend only on the eigenvalues of $J_f$ and thus the above equation gives a relation on the spectra of the differentials of $f$ at the fixed points. On the other hand, the Lefschetz number $\operatorname{L}(f,\wedge^k df,\Omega^k_X)$ depends only on the action of $f$ on the Dolbeault cohomology of $X$ in virtue of the standard identification $H^q(X,\Omega^p_X)\cong H^{p,q}_{\bar{\partial}}(X)$. In the particular case $X=\P^n$ this action on the Dolbeault cohomology depends only on the degree $d$ of $f$, and in fact, it is not hard to see that it is just multiplication by an appropiate power of $d$. This gives the following result.

\begin{corollary}\label{coro:indicesPn}
Let $f$ be a transversal holomorphic endomorphism of $\P^n$ of degree $d$. Then
\begin{equation}\label{eq:LefschetzPn}
\sum_{f(p)=p}\frac{\tr(\wedge^k J_f(p))}{\det(I-J_f(p))} = (-1)^k d^k .
\end{equation}
\end{corollary}

\begin{remark}
The Woods Hole formula can also be used to count the number of fixed points or singular points in different settings. For example, we can choose $\Fs=\wedge^\ast\Omega^1_{\P^n}$, the whole exterior algebra of $\Omega^1_{\P^n}$ with $\varphi$ acting as $\wedge^k df$ on the graded componen of degree $k$. In such case, it follows from \hyperref[coro:indicesPn]{Corollary \ref*{coro:indicesPn}} and from the fact that $\det(I-A) = \sum_k (-1)^k \tr(\wedge^k A)$ that a degree $d$ endomorphism of $\P^n$ has $1+d+\ldots+d^n=\frac{d^{n+1}-1}{d-1}$ fixed points.
\end{remark}

\section{Indices for transversal endomorphisms of \texorpdfstring{$\P^n$}{Pn}}\label{sec:endomorphisms}

Besides the index theorems stated in \hyperref[coro:indicesPn]{Corollary \ref*{coro:indicesPn}} we can use the Woods Hole formula to prove the following index theorem.

Recall that any invariant polynomial function on $\mathfrak{gl}_n(\C)$ can be uniquely expressed as a polynomial on the \emph{elementary symmetric polynomials} $\sigma_1, \ldots, \sigma_n$.

\begin{application}\label{app:indicesPn}
Let $B$ be an invariant polynomial function on $\mathfrak{gl}_n(\C)$ of degree at most $n$. Let $Q\in\C[z_1,\ldots,z_n]$ be the unique polynomail that satisfies
\[ B=Q(\sigma_1,\ldots,\sigma_n). \]
Then if $f$ is a transversal endomorphism of $\P^n$ of degree $d$ we have
\begin{equation}\label{eq:indicesPn} 
\sum_{f(p)=p}\frac{B(J_f(p))}{\det(I-J_f(p))} = Q(-d,(-d)^2,\ldots,(-d)^n).
\end{equation}
\end{application}

Note that if in particular $B=\sigma_k$ then the right hand side of the above formula is simply $(-1)^k d^k$ and we recover equation (\ref{eq:LefschetzPn}). This result refines the main result in \cite{Guillot2004} by computing the value of the sum of indicies, and is also in accordance with the results in \cite{Ueda1995,Abate2014}.

\begin{proof}[Sketch of the proof]
Both side of equation (\ref{eq:indicesPn}) depend linearly on the polynomial $B$, hence it is enough to prove the equality in the case $B$ is a monic monomial. Let us assume then that $B$ is given by
\[ B=\sigma_1^{a_1}\ldots\sigma_n^{a_n}. \]
Since $B$ has degree at most $n$ we must have $a_1+2a_2+\ldots+na_n\leq n$. Define \[ \Fs=(\Omega^1_{\P^n})^{\otimes a_1}\otimes\ldots\otimes(\Omega^n_{\P^n})^{\otimes a_n}, \]
and $\varphi=(df)^{\otimes a_1}\otimes\ldots\otimes(\wedge^n df)^{\otimes a_n} \colon f^*\Fs\to\Fs$ and apply the Woods Hole formula. Since the trace of a tensor product is the product of the traces it is clear that the trace of $\varphi(p)$ acting on $\Fs(p)$ will be precisely $B(J_f(p))$. Let $k=\operatorname{deg}\,B=a_1+2a_2+\ldots+na_n$. We need to prove now that the Lefschetz number $\operatorname{L}(f,\varphi,\Fs)$ equals $Q(-d,(-d)^2,\ldots,(-d)^n)=(-d)^{a_1}\ldots(-d)^{na_n}=(-d)^k$. Formula (\ref{eq:indicesPn}) follows immediately from the next lemma.

\begin{lemma}\label{lemma:cohomology}
If $k=a_1+2a_2+\ldots+na_n\leq n$ then the natural map 
\[ \wedge \colon \Fs=(\Omega^1_{\P^n})^{\otimes a_1}\otimes\ldots\otimes(\Omega^n_{\P^n})^{\otimes a_n} \longrightarrow \Omega^k_{\P^n} \]
induces isomorphisms in all cohomology groups. In particular, the cohomology of $\Fs$ vanishes in all degrees except $k$ and $H^k(\P^n,\Fs)\cong\C$. Moreover the following diagram commutes:
\[ \xymatrix{
H^k(\P^n,\Fs) \ar[r]^{\varphi} \ar[d]_{\wedge} & H^k(\P^n,\Fs) \ar[d]^{\wedge} \\
H^k(\P^n,\Omega^k_{\P^n}) \ar[r]_{\wedge^k df}  & H^k(\P^n,\Omega^k_{\P^n}) .
 }
\]
\end{lemma}

Since $\wedge\colon H^k(\P^n,\Fs)\to H^k(\P^n,\Omega^k_{\P^n})$ is an isomorphism it follows from the commutativity of the above diagram that $\tr(\varphi,\,H^k(\P^n,\Fs))=\tr(\wedge^k df,\,H^k(\P^n,\Omega^k_{\P^n}))=d^k$ and so $\operatorname{L}(f,\varphi,\Fs)=(-1)^k d^k$, which is what we wanted to prove.
\end{proof}

The main step in the proof of \hyperref[lemma:cohomology]{Lemma \ref*{lemma:cohomology}} is proving that $\Fs$ actually has the same cohomology groups as $\Omega^k_{\P^n}$ (once this is achieved it is not hard to see that the map $\wedge$ induces an isomorphism in degree $k$ cohomology which is the only non-trivial degree). Such step is proved by an induction argument -- we shall prove it here only for $\Fs=\Omega^p_{\P^n}\otimes\Omega^q_{\P^n}$.

\begin{proof}[Proof of Lemma \ref*{lemma:cohomology}]
Let us prove that $\Omega^p_{\P^n}\otimes\Omega^q_{\P^n}$ and $\Omega^{p+q}_{\P^n}$ have isomorphic cohomology groups by inducting on $p$. The statement is clearly true for $p=0$. Assume now that it holds for $p-1$ and consider the generalized Euler sequence 
\begin{equation}\label{eq:Eulersequence}
0 \longrightarrow \Omega^p_{\P^n} \longrightarrow E\otimes \O_{\P^n}(-p) \longrightarrow \Omega^{p-1}_{\P^n} \longrightarrow 0,
\end{equation}
where $E$ denotes the trivial vector bundle of rank $\binom{n+1}{p}$. Since $\Omega^q_{\P^n}$ is locally free, the stalks $(\Omega^q_{\P^n})_x$ are flat $\O_{\P^n,x}$-modules and so tensoring (\ref{eq:Eulersequence}) with the sheaf $\Omega^q_{\P^n}$ preserves the exactness of the sequence. We obtain
\[ 0 \longrightarrow \Omega^p_{\P^n}\otimes\Omega^q_{\P^n} \longrightarrow E\otimes \Omega^q_{\P^n}(-p) \longrightarrow \Omega^{p-1}_{\P^n}\otimes\Omega^q_{\P^n} \longrightarrow 0. \]
It is known that if $p+q\leq n$ then $\Omega^q_{\P^n}(-p)$, and thus $E\otimes\Omega^q_{\P^n}(-p)$, has vanishing cohomology in every dimension. This implies that the connecting homomorphisms 
\[ \delta\colon H^k(\P^n,\Omega^{p-1}_{\P^n}\otimes\Omega^q_{\P^n})\to H^{k+1}(\P^n,\Omega^p_{\P^n}\otimes\Omega^q_{\P^n}) \]
are isomorphisms for all $k$. This proves that that $\Omega^p_{\P^n}\otimes\Omega^q_{\P^n}$ and $\Omega^{p+q}_{\P^n}$ have isomorphic cohomology groups.

Lastly, the composition 
\[ \xymatrix{
H^p(\P^n,\Omega^p_{\P^n}) \otimes H^q(\P^n,\Omega^q_{\P^n}) \ar[r]^-{\cup} & H^{p+q}(\P^n,\Omega^p_{\P^n}\otimes \Omega^q_{\P^n}) \ar[r]^-{\wedge} & H^{p+q}(\P^n,\Omega^{p+q}_{\P^n})
 }
\]
of the cup product and the homomorphism induced by $\wedge$ can be understood as the intersection of analytic cycles in the Chow ring of projective space and is thus an isomorphism. We conclude that $\wedge$ induces isomorphisms in $H^{p+q}$ and hence in all cohomology groups.
\end{proof}

\section{The Baum-Bott theorem for foliations on \texorpdfstring{$\P^2$}{P2}} \label{sec:foliations}

There are two important index theorems for foliations on $\P^2$: the Baum-Bott theorem and the Camacho-Sad theorem for invariant curves. The former can be reduced, under a non-degeneracy assumption, to a particular case of the Woods Hole formula. In order to do so we need to be able to associate to any holomorphic foliation an endomorphism of $\P^2$.

\begin{remark}
Even though we will be able to associate an endomorphism of $\P^2$ to any holomorphic foliation, the choice will not be unique. Also, it is important to point out that a curve that is invariant for the foliation need not be invariant for the associated endomorphism. Indeed, invariant curves are extremely rare for endomorphisms and in particular a smooth, totally invariant curve has to be a line \cite{FornaessSibony1994,CerveauLinsNeto2000}. In the very particular case of an exceptional line we can recover the Camacho-Sad theorem from the Woods Hole formula.
\end{remark}

We now proceed as follows: let $\F$ be a holomorphic foliation on $\P^2$ having only isolated non-degenerate singularities. We can always choose a homogeneous vector field $\xi=P_0\frac{\partial}{\partial x_0} + P_1\frac{\partial}{\partial x_1} + P_2\frac{\partial}{\partial x_2}$ on $\C^3\setminus\{0\}$ that descends to a (singular) line field on $\P^2$ which induces foliation $\F$. The choice is not unique, but every two such choices differ only by a multiple of the radial vector field. Define an endomorphism $f_\xi$ of $\P^2$ by setting
\[ f_\xi(x_0,x_1,x_2) = [P_0:P_1:P_2] .\]
It is not hard to see that the fixed points of $f_\xi$ correspond to the singularities of $\F$ and thus all these are simple and isolated. Moreover we have the following fundamental observations:

\begin{proposition}\label{prop:fundamental}
In affine coordinates $(x,y)=(\frac{x_1}{x_0},\frac{x_2}{x_0})$ foliation $\F$ is induced by
\[ (P_1(1,x,y)-xP_0(1,x,y)\frac{\partial}{\partial x} + (P_2(1,x,y)-yP_0(1,x,y)\frac{\partial}{\partial y}. \]
Therefore if $p=[p_0:p_1:p_2]$ is a singular point of $\F$ with $p_0\neq 0$ there exists a small open neighborhood of $p$ where $P_0(1,x,y)$ does not vanish and so foliation $\F$ is induced by the vector field
\[ \tilde{\xi} = \left(\frac{P_1(1,x,y)}{P_0(1,x,y)}-x\right)\frac{\partial}{\partial x} + \left(\frac{P_2(1,x,y)}{P_0(1,x,y)}-y\right)\frac{\partial}{\partial y}, \]
locally around $p$.
\end{proposition}

Note that around such a singular point $p$ the endomorphism $f_\xi$ is locally given by $(x,y)\mapsto (\frac{P_1}{P_0},\frac{P_2}{P_0})$ and thus we have the relation
\[ D\tilde{\xi}(p)=Df_\xi (p)-I. \]

\begin{application}[Baum-Bott theorem]\label{app:BB}
Let $\F$ be a homlomorphic foliation of degree $d$ on $\P^2$ having isolated non-degenerate singularities only. Then
\[ \sum_{\operatorname{Sing }\F}\operatorname{BB}(\F,p)=(d+2)^2. \]
\end{application}

\begin{proof}
Let $\Fs=\Omega^1_{\P^2}\otimes\Omega^1_{\P^2}$ and $\varphi$ the lift of $f_\xi$ to $\Fs$ given by $\varphi=df_\xi\otimes df_\xi$ and apply the Woods Hole formula.
\end{proof}

\section{Indices for polynomial vector fields on \texorpdfstring{$\C^2$}{C2}} \label{sec:vectorfields}

We have seen in the previous section that whenever we have a foliation $\F$ on $\P^2$ we can always associate to it a homogeneous vector field $\xi$ on $\C^3$ and an endomorphism $f_\xi$ of $\P^2$, however, this endomorphism is not uniquely determined by $\F$. Suppose now that we are given a degree $d$ polynomial vector field $v=P(x,y)\frac{\partial}{\partial x}+Q(x,y)\frac{\partial}{\partial y}$ on $\C^2$. In this case there is a natural choice for the associated endomorphism given by
\begin{equation} \label{eq:fv}
 f_v(x_0,x_1,x_2) = [x_0^d:\tilde{P}(x_0,x_1,x_2)+x_0^{d-1}x_1:\tilde{Q}(x_0,x_1,x_2)+x_0^{d-1}x_2], 
\end{equation}
where $\tilde{P}$ and $\tilde{Q}$ are the homogenizations of $P,Q$ with respect to the variable $x_0$. At every point $p$ on $\C^2$ we have the relation
\[ Dv(p) = Df_v(p)-I. \]
Moreover, if $v$ has an invariant line at infinity and non-degenerate singularities then $f_v$ has simple isolated fixed points only. With this construction we can now prove the Euler-Jacobi formulae as a particular case of the Woods Hole formula.

\begin{application}[$1^{\mathrm{st}}$ Euler-Jacobi relation]\label{app:EJ1}
Let $v$ be a polynomial foliation on $\C^2$ of degree $d$ having $d^2$ non-degenerate singularities. Then
\[ \sum_{v(p)=0}\frac{1}{\det(Dv(p))}=0. \]
\end{application}

\begin{proof}
Let $L$ denote the line at infinity and let $\mathcal{I}_L$ denote the sheaf of holomorphic functions on $\P^2$ that vanish on $L$. Under the hypotheses of the theorem the endomorphism $f_v$ has degree $d$, preserves the line $L$, and the sheaf $f^*\mathcal{I}_L$ equals $\mathcal{I}_L^{\otimes d}$ (i.e.~the sheaf of functions of vanishing order at least $d$ along $L$). There is a natural inclusion $\iota\colon f^*\mathcal{I}_L\cong \mathcal{I}_L^{\otimes d} \to \mathcal{I}_L$. Applying the Woods Hole formula to the particular case $\Fs=\mathcal{I}_L$, $\varphi=\iota$ gives the Euler-Jacobi formula stated above.
\end{proof}

\begin{application}[$2^{\mathrm{nd}}$ Euler-Jacobi relation]\label{app:EJ2}
Let $v$ be a polynomial foliation on $\C^2$ of degree $d$ having $d^2$ non-degenerate singularities. Then
\[ \sum_{v(p)=0}\frac{\tr(Dv(p))}{\det(Dv(p))}=0. \]
\end{application}

\begin{proof}
Consider $\Fs=\Omega^1_{\P^2}\otimes\mathcal{I}_L$, the sheaf of 1-forms on $\P^2$ that are identically zero on $L$. The second Euler-Jacobi relation follows from the Woods Hole formula in the particular case where $\Fs=\Omega^1_{\P^2}\otimes\mathcal{I}_L$, $\varphi = df\otimes\iota$.
\end{proof}

Again, let $L$ denote the line at infinity and assume $v$ is non-dicritic. From the expression for $f_v$ in (\ref{eq:fv}) it is clear that the line $L$ is totally invariant for the endomorphism $f_v$ (that is, $f_v^{-1}(L)=L$). Under this special situation we can recover the Camacho-Sad theorem from the Woods Hole formula.

\begin{application}[Camacho-Sad theorem, particular case]\label{app:CS}
If the foliation $\F_v$ defined by $v$ on $\P^2$ has non-degenerate singularities along the invariant line at infinity $L$ we have
\[ \sum_{L\cap\operatorname{Sing }\F_v}\operatorname{CS}(\F_v,L,p)=1. \]
\end{application}

\begin{proof}
Let $N_L$ denote the normal bundle of $L$ in $\P^2$ and recall that $\operatorname{deg }N_L=N_L\cdot N_L=1$. Let $\mathcal{N}^*_L$ denote the conormal sheaf of $L$, that is, the sheaf of sections of the dual of the normal bundle $N_L$. A section of $\mathcal{N}^*_L$ can be understood as a holomorphic 1-form on $\P^2$ that vanishes along the tangent bundle of $L$. If $L$ is totally invariant by $\F$ we have that $f_\xi(L)=L$ and if $\omega$ is a 1-form that vanishes along $TL$ so does $f^*\omega$. This implies that $df_\xi$ acts on $\mathcal{N}^*_L$ in a natural way. Let us apply Woods Hole on the line $L$ by making $\Fs=\mathcal{N}^*_L$ and lifting $f_\xi|_L$ to $\mathcal{N}^*_L$ by $\varphi=df_\xi$. The Woods Hole formula applied to this particular case implies the Camacho-Sad theorem.
\end{proof}

\section*{Table: Summary of applications}

\renewcommand{\arraystretch}{1.25}
\begin{center}
 \begin{tabular}{ | l | c | c | c | c |} 
 \hline
Formula & Manifold $X$ & Map $f$ & Sheaf $\Fs$ & Lift $\varphi$ \\ \hline
 Generalized Lefschetz & Any $X$ & $f$ transversal & $\Omega^k_X$ & $\wedge^k df$ \\ \hline
 Guillot's relations & $\P^n$ & $f$ transversal & $\bigotimes(\Omega^k_{\P^n})^{\otimes a_k}$ & $\bigotimes (\wedge^k df)^{\otimes a_k}$ \\ \hline
 Baum-Bott & $\P^2$ & $f_\xi$ & $\Omega^1_{\P^2}\otimes\Omega^1_{\P^2}$ & $df_\xi\otimes df_\xi$ \\ \hline
 Euler-Jacobi 1 & $\P^2$ & $f_v$ & $\mathcal{I}_L$ & $\iota$ \\ \hline
 Euler-Jacobi 2 & $\P^2$ & $f_v$ & $\Omega^1_{\P^2}\otimes\mathcal{I}_L$ & $df\otimes\iota$ \\ \hline
 Camacho-Sad & $L\hookrightarrow\P^2$ & $f_\xi|_L$ & $\mathcal{N}^*_L$ & $df_\xi$ \\ \hline
 \end{tabular}
\end{center}

\vspace{4mm}

\bibliographystyle{alpha}
\bibliography{ref_WH}

\begin{thebibliography}{CLN00}

\bibitem[AB66]{AtiyahBott1966}
M.~F. Atiyah and R.~Bott.
\newblock A {L}efschetz fixed point formula for elliptic differential
  operators.
\newblock {\em Bull. Amer. Math. Soc.}, 72:245--250, 1966.
\newblock
  \href{http://dx.doi.org/10.1090/S0002-9904-1966-11483-0}{DOI:10.1090/S0002-9904-1966-11483-0}.

\bibitem[AB67]{AtiyahBott1967}
M.~F. Atiyah and R.~Bott.
\newblock A {L}efschetz fixed point formula for elliptic complexes. {I}.
\newblock {\em Ann. of Math. (2)}, 86:374--407, 1967.
\newblock \href{http://dx.doi.org/10.2307/1970694}{DOI:10.2307/1970694}.

\bibitem[AB68]{AtiyahBott1968}
M.~F. Atiyah and R.~Bott.
\newblock A {L}efschetz fixed point formula for elliptic complexes. {II}.
  {A}pplications.
\newblock {\em Ann. of Math. (2)}, 88:451--491, 1968.
\newblock \href{http://dx.doi.org/10.2307/1970721}{DOI:10.2307/1970721}.

\bibitem[Aba14]{Abate2014}
M.~Abate.
\newblock Index theorems for meromorphic self-maps of the projective space.
\newblock In {\em Frontiers in complex dynamics}, volume~51 of {\em Princeton
  Math. Ser.}, pages 451--461. Princeton Univ. Press, Princeton, NJ, 2014.
\newblock \href{http://www.ams.org/mathscinet-getitem?mr=3289918}{MR:3289918}.
  Preprint available at \href{http://arxiv.org/abs/1106.2394}{arXiv:1106.2394}.

\bibitem[Bea72]{Beauville1972}
A.~Beauville.
\newblock {\em Formules de points fixes en cohomologie coherente}.
\newblock Seminaire de Geometrie Algebraique. Secrétariat mathématique de
  l'ÉNS, 1972.
\newblock Available online at
  \href{http://math1.unice.fr/~beauvill/pubs/lef.pdf}{http://math1.unice.fr/~beauvill/pubs/lef.pdf}.

\bibitem[CLN00]{CerveauLinsNeto2000}
D.~Cerveau and A.~Lins~Neto.
\newblock Hypersurfaces exceptionnelles des endomorphismes de {${\bf C}{\rm
  P}(n)$}.
\newblock {\em Bol. Soc. Brasil. Mat. (N.S.)}, 31(2):155--161, 2000.
\newblock \href{http://dx.doi.org/10.1007/BF01244241}{DOI:10.1007/BF01244241}.

\bibitem[FS94]{FornaessSibony1994}
J.~E. Forn{\ae}ss and N.~Sibony.
\newblock Complex dynamics in higher dimension. {I}.
\newblock {\em Ast\'erisque}, (222):5, 201--231, 1994.
\newblock \href{http://www.ams.org/mathscinet-getitem?mr=1285389}{MR:1285389}.

\bibitem[Gui04]{Guillot2004}
A.~Guillot.
\newblock Un th\'eor\`eme de point fixe pour les endomorphismes de l'espace
  projectif avec des applications aux feuilletages alg\'ebriques.
\newblock {\em Bull. Braz. Math. Soc. (N.S.)}, 35(3):345--362, 2004.
\newblock
  \href{http://dx.doi.org/10.1007/s00574-004-0018-7}{DOI:10.1007/s00574-004-0018-7}.

\bibitem[Ram16]{TwinVectorFields}
V.~Ram\'{i}rez.
\newblock Twin vector fields and independence of spectra for quadratic vector
  fields.
\newblock {\em \normalfont{To appear in:} \textit{J. Dynam. Control Systems}},
  2016.
\newblock Preprint available at
  \href{http://arxiv.org/abs/1508.02413}{arXiv:1508.02413v3}.

\bibitem[Tae15]{Taelman2015}
L.~Taelman.
\newblock {\em Sheaves and Functions Modulo $p$}.
\newblock London Mathematical Society Lecture Note Series (No.~429). Cambridge
  University Press, 2015.
\newblock
  \href{http://dx.doi.org/10.1017/CBO9781316480687}{DOI:10.1017/CBO9781316480687}.

\bibitem[Ued95]{Ueda1995}
T.~Ueda.
\newblock Complex dynamics on projective spaces---index formula for fixed
  points.
\newblock In {\em Dynamical systems and chaos, {V}ol. 1 ({H}achioji, 1994)},
  pages 252--259. World Sci. Publ., River Edge, NJ, 1995.
\newblock \href{http://www.ams.org/mathscinet-getitem?mr=1479941}{MR:1479941}.

\bibitem[Ver67]{Verdier1967}
J.~L. Verdier.
\newblock The {L}efschetz fixed point formula in etale cohomology.
\newblock In {\em Proc. {C}onf. {L}ocal {F}ields ({D}riebergen, 1966)}, pages
  199--214. Springer, Berlin, 1967.
\newblock \href{http://www.ams.org/mathscinet-getitem?mr=0242846}{MR:0242846}.

\end{thebibliography}

\end{document}